


\magnification=\magstep 1
\parindent = 0 pt
\baselineskip = 16 pt
\parskip = \the\baselineskip

\font\AMSBoldBlackboard = msbm10

\def\RR{{\hbox{\AMSBoldBlackboard R}}}
\def\CC{{\hbox{\AMSBoldBlackboard C}}}
\def\AA{{\hbox{\AMSBoldBlackboard A}}}
\def\QQ{{\hbox{\AMSBoldBlackboard Q}}}

\def\ZZ{{\hbox{\AMSBoldBlackboard Z}}}

\settabs 12\columns
\rightline{math.NT/9810169}
\vskip 1 true in
{\bf \centerline{THE EXPLICIT FORMULA IN SIMPLE TERMS}}
\vskip 0.5 true in
\centerline{Jean-Fran\c{c}ois Burnol}
\par
\centerline{October 1998}
\centerline{revised November 1998}
\par
This is a semi-expository paper on the easier aspects of the Explicit Formula for the Riemann Zeta Function. The topics reviewed here include: Weil's criterion for the Riemann Hypothesis and its probabilistic interpretation, various formulations of the contribution corresponding to the real place, Haran's version of the Explicit Formula, and the author's own derivation which puts all places on the same footing. This derivation, an addendum to Tate's Thesis, is in the spirit of Weil's insights towards an adelic understanding of the Explicit Formula. Whereas the analyst would likely formulate the Explicit Formula as a multiplicative convolution, Haran's theorem shows it is also an additive convolution. An intriguing conductor operator is considered whose spectral analysis is equivalent to the Explicit Formula. At a finite place it has a positive cuspidal spectrum. I also comment on the flexibility still left in the Explicit Formula, which shows that it has a symmetry group containing the non-zero rational numbers.
\vfill
{\parskip = 0 pt
62 rue Albert Joly\par
F-78000 Versailles\par
France\par}

jf.burnol@dial.oleane.com

\eject

{\bf
TABLE OF CONTENTS\par
\par
The formula of Riemann\thinspace-\thinspace von Mangoldt as a distribution\par
Weil's positivity criterion and stochastic processes\par
The explicit formula: classical approach\par
A convolution algebra\par
Weil's Theorem\par
The functional equation and homogeneous distributions\par
Haran's Theorem\par
Mellin and Fourier\par
The explicit formula: adelic approach\par
The conductor operator\par
The symmetries of the explicit formula\par
References\par}
\vfill \eject

{\bf The formula of Riemann\thinspace-\thinspace von Mangoldt as a distribution}\par

The paper of Riemann [1] unveiled a deep relationship between the prime numbers and the zeros in the complex plane of a function already considered by Euler, the famous zeta function. A striking explicit formula was given by Riemann for the prime-counting function in terms of the zeros, and the generalizations obtained later all go by the name of ``Explicit Formulae". This classical theory is surveyed in the books of Edwards [2] and Ingham [3]. A particularly elegant formula was rigorously proven by von Mangoldt:
$$\sum_{1<n<X}{\Lambda(n)} + {1\over2}\Lambda(X) = 
X - \sum_{\rho}{X^\rho \over\rho} - \log(2\pi) - {1\over2}\log(1-X^{-2})$$
Here $X > 1$ (not necessarily an integer) and $\Lambda(Y) = \log(p)$ if $Y > 1$ is a positive power of the prime number $p$, and is $0$ for all other values of $Y$. The $\rho$'s are the Riemann Zeros (in the critical strip), the sum over them is not absolutely convergent, even after pairing $\rho$ with $1-\rho$. It is defined as $\lim_{T\to\infty} \sum_{|Im(\rho)| < T}{X^\rho \over\rho}$.

One can not directly plug $X = 1$ in this equation, but nevertheless one has another explicit evaluation, this time for an absolutely convergent sum:
$$\sum_{\rho}{1\over{\rho(1-\rho)}} = 2 + \gamma - \log(4\pi)$$
which follows (for example) from the functional equation combined with the Hadamard product representation for the zeta function (see [2] chap. 3, or Patterson [4], chap. 3).

I can not help mentioning here a little joke reformulation of the Riemann Hypothesis.

{\bf Lemma:} The Riemann Hypothesis holds if and only if the following evaluation holds:
$$\sum_{\rho}{1\over{|\rho|^2}} = 2 + \gamma - \log(4\pi)$$

{\bf Proof:} This statement follows from the equality:
$${1\over2} ({1\over|\rho|^2} + {1\over|1-\rho|^2}) = Re({1\over\rho(1-\rho)}) + {1\over2} {{(1-2Re(\rho))}^2\over|\rho|^2|1-\rho|^2}$$

The literature on the Explicit Formulae is vast, it seems (to the best of the author's knowledge) that Weil's version from 1952 [5] subsumes most previously known such formulae. Its true merit probably lies in the adelic direction as will be explained later in this paper, although it is often advertised as having introduced a characteristic distribution-theoretic flavor to the subject. It also introduced a funny psychological change in the understanding of what ``explicit" stands for. In the Riemann-van Mangoldt result one expresses an arithmetical counting function in an explicit manner in terms of the zeros. This is why the primes stand on the left and the zeros on the right of the equality sign. In Weil's tradition one expresses a sum over the zeros in an explicit manner in terms of the primes, hence the zeros are on the left and the primes on the right! Rieman wanted to say things about the primes from the zeros (and indeed this later led to a proof of the Prime Number Theorem), and Weil wanted to say things about the zeros from (some construction based on) the primes.

Let $g$ be a function on the positive half-line. Its Mellin transform is
$$\widehat{g}(s) = \int_{0}^{\infty} g(u){u}^s \, { du \over u}$$
and Weil's formula is
$$\widehat{g}(0) + \widehat{g}(1) - \sum_{\rho}{\widehat{g}(\rho)}
= \sum_\nu{W_\nu(g)}$$
Here $\nu$ indexes the valuations of the number field $\QQ$ (the primes $p$ and the archimedean place $r$ corresponding to the real numbers). The $W_\nu$'s are distributions on the positive half-line, whose explicit forms will be given later. For a finite prime $p$, $W_p$ is supported on powers and inverse powers of $p$ (except $1$), whereas $W_r$ has a singularity at $1$ but is otherwise smooth. The singularity at $1$ would reappear even for a finite place if we considered more general Dirichlet L-series instead of the Riemann Zeta function, so that we can imagine that the real place always has some ramification.

Of course one has to be specific about the functions $g$ allowed, and for simplicity sake we could restrict to the space $\cal D$ of ``test functions'' $g$ (smooth and compactly supported), in which case all the involved sums and integrals are absolutely convergent (one has $|\widehat{g}(s)| = O((1 + |s|)^{-N})$ for any fixed $N$, uniformly in each vertical strip $A\leq Re(s)\leq B$). But Weil proved his formula for a large class of $g$'s, including in particular the step-function for which $\widehat{g}(s) = {(X^s - 1)\over s}$, so that it is a genuine generalization of von\thinspace Mangoldt's theorem and requires in general the subtle definition mentioned above for the sum over the zeros (it is an interesting exercise to actually exhibit von\thinspace Mangoldt's result as a special case because of the delicate expression for $W_r(g)$ given by Weil).

In fact we could imagine trying to prove Weil's formula by first generalizing von\thinspace Mangoldt's formula to all step functions, then applying a suitable density argument, but this quickly leads into complicated matters because of the subtle nature of the convergence. Note that this is a subtlety linked with the inversion problem for Fourier series or integrals and has nothing to do with primes {\it per se}.

The actual location of the zeros is well hidden in Weil's Formula. Its derivation uses only the basic properties of the zeta function and of the gamma function necessary for the proof of the Hadamard Product Formula. It is thus less deep than the Prime Number Theorem as it does not require the non vanishing on $Re(s) = 1$. As the joke reformulation of the Riemann Hypothesis shows, any explicit formula for a {\it non-holomorphic} summand would presumably immediately tell us much more on the zeros. But Weil has found a neat way to express the Riemann Hypothesis which greatly motivates the search for a better understanding of the $W_\nu$'s. This is the topic of the next section.

{\bf Weil's positivity criterion and stochastic processes\par}

{\bf Lemma:} The Riemann Hypothesis holds if and only if 
$$\sum_{\rho}{\widehat{g}(\rho)\cdot\overline{\widehat{g}(\overline{1-\rho})}} \geq 0$$
for all smooth compactly supported functions $g$ on $(0, \infty)$.

{\bf Proof:} As $\rho =\overline{1-\rho}$ on the critical line one direction is clearly true. For the other one, Weil's proof uses a larger class of functions $g$ (he needs $\widehat{g}(s) = P(s){e^{-A(s-1/2)^2}}$ for polynomials $P(s)$). Rather than applying a suitable density argument, it is easier to prove the Lemma directly as stated.

First, there is the easy fact that given $N$ distinct complex numbers $z_1, z_2, ...., z_N$ and $N$ values $a_1, a_2, ..., a_N$ there is $g \in \cal D$ such that $\widehat{g}(z_i) = a_i$. If this were false there would be a linear relation among the linear forms $g \mapsto \widehat{g}(z_i)$. As $\widehat{Dg}(s) = s\widehat{g}(s)$ where $(Dg)(x) = - xg'(x)$, we can eliminate one of the $z_i$'s and get a relation involving $N - 1$ linear forms. One ends up with $N = 1$ where the situation is clear.

So, let $\rho_0$ be a bad zero, so that $\rho_0 \neq \overline{1-\rho_0}$. We pick $g$ such that $\widehat{g}(\rho_0) = \widehat{g}(\overline{1-\rho_0}) = 1$. There will be only finitely many other (``non-trivial'') zeros where $|\widehat{g}(\rho)| \geq 1/2$. Let's take now $k$ with $\widehat{k}$ vanishing at all those zeros and taking the values $1$ at $\rho_0$ and $-1$ at $\overline{1-\rho_0}$. Defining $g_N \in \cal D$ to be the multiplicative convolution of $g$ ($N$ times) and of $k$, we find
$$\sum_{\rho}{\widehat{g_N}(\rho)\cdot\overline{\widehat{g_N}(\overline{1-\rho})}} = -2 + O(1/4^N)$$
so that it is negative for $N$ large enough. This completes the proof of Weil's positivity criterion.

{\bf Note:} We can choose $k$ so that $\widehat{k}(0) = \widehat{k}(1) = 0$ so that we also have $\widehat{g_N}(0) = \widehat{g_N}(1) = 0$.

As was alluded to before the main impulse about this positivity criterion has been to look for a suitable {\it geometric} understanding. This is natural if one takes into account the well-known function-field story. But another road can be contemplated: that of {\it randomness}. For this one first shifts the critical line by $1/2$, and moves the contribution of the poles to the same side of the explicit formula as the primes. So the explicit formula is now seen (as in Weil's paper [5]) to give a distribution $C$ on $(0, \infty)$ such that
$$C(g) = \sum_{\rho}{\widehat{g}(\rho - {1\over2})}$$
and Weil's positivity criterion becomes
$$C(g*g^c) \geq 0$$
where $g^c = \overline{g({1\over u})}$.

One obvious way ([6]) to obtain this positivity is to imagine that we have a probability measure $\mu$ on the (tempered) distributions $\Omega$ on $(0, \infty)$ with the $\sigma-$algebra generated by the cylinder sets, invariant under multiplicative translations, and the identity:
$$\int_{\Omega}X_F(d)\overline{X_G(d)}\, d\mu(d) = C(F*G^c)$$
where $F$ and $G$ are two arbitrary elements of $\cal D$ and $X_F$ and $X_G$ are the two associated ``coordinates'' on the space of distributions $d \mapsto X_F(d) = d(F)$, $d \mapsto X_G(d) = d(G)$. For $F = G$ it is obviously non-negative!\par

We also assume that $\int_{\Omega}X_F(d)\, d\mu(d) = 0$ for all ``coordinates'' $X_F$, so that $C(F*G^c)$ is the covariance of the random variables defined by $F$ and $G$. It is known (see [7]) that the covariance of any generalized stochastic process invariant under time translation and with zero mean has the shape $C(F*G^c)$ for some distribution $C$ (in general of course this is considered on the additive real line, so that $*$ is an additive convolution and $G^c(x) = \overline{G(-x)}$\thinspace).

The general case of Dirichlet L-series leads to expect an all-encompassing generalized stochastic process with ``time'' given by the idele classes of $\QQ$ (which is where Weil's distribution $C$ naturally lives.) 

If such an ``arithmetic stochastic process'' could be constructed, then the Riemann Hypothesis for all Dirichlet L-series would follow. We could hope to identify this random process through its characteristic functional $E(F) = \int_{\Omega}\exp(iX_F(d))\, d\mu(d)$, perhaps by comparing with the function-field case, but I haven't been able to go very far in this direction, so that I can not really advocate further this probabilistic point of view here. There are misgivings about the necessary shift by 1/2 and the fact of putting the contributions of the poles and of the primes together; they arise from a later topic of this paper, the adelic formulation of the Explicit Formula.

{\bf The explicit formula: classical approach\par}

Let $K_\varepsilon (u) = u^{-\varepsilon} + u^{1 + \varepsilon}$ for some $\varepsilon > 0$. Let $g$ be a Lebesgue measurable function on $(0, \infty)$ such that condition $(A_\varepsilon)$ holds:
$$g(u)K_\varepsilon (u) \hbox{ has bounded total variation on }(0, \infty)\leqno (A_\varepsilon)$$

{\bf Lemma:} Under condition $(A_\varepsilon)$\ $\widehat{g}(s)$ exists as an analytic function in the open strip $-\varepsilon < Re(s) < 1+\varepsilon$, and satisfies $|\widehat{g}(s)| = O(1/|Im(s)|)$ uniformly in this strip.

{\bf Proof:} A well-known lemma (see [4], app. 1) tells us that for $f(t)$ in $L^1(\RR, dt)$, vanishing at infinity, $|\int_{-\infty}^{\infty} f(t)e^{iTt}\, dt| \leq {1\over T} V(f)$ with $V(f)$ the total variation of $f$. So $|\widehat{g}(c+iT)| \leq {1\over T} V_0^\infty(g(u)u^c)$ for $-\varepsilon < c < 1+\varepsilon$. Now another lemma $V(fg) \leq \sup(|g|)V(f) + \sup(|f|)V(g)$ implies $$V_0^\infty(g(u)u^c) \leq \sup(|g(u)K_\varepsilon (u)|)V_0^\infty(u^c/K_\varepsilon (u)) + V_0^\infty(g(u)K_\varepsilon (u))$$
and this is $O(1)$ as $V_0^\infty(u^c/K_\varepsilon (u))\leq 2$.

Furthermore $g$ itself will have bounded total variation so that the right-limit $g(u+)$ and left-limit $g(u-)$ exist and we also assume condition $(B)$:
$$\forall u\quad g(u+) + g(u-) = 2g(u)\leqno (B)$$

{\bf Theorem:} Under conditions $(A_\varepsilon)$ and $(B)$
$$\lim_{T\to\infty} {1\over{2\pi i}}\int_{c-iT}^{c+iT}-{\zeta'(s)\over\zeta(s)}\widehat{g}(s)\, ds = \sum_p V_p(g)$$
for any $c$ with $1 < c < 1+\varepsilon$, with $V_p(g) := \log(p)\sum_{k \geq 1} g(p^k)$ (absolutely convergent).

{\bf Proof:} Note that the integral is not necessarily absolutely convergent. The trick (from Weil [5]) is to work in reverse, defining first $G(u) = \sum_{p, k\geq1}\log(p) \,g(p^k u)$. One checks in succession:
\+\cr
\+&The sum defining $G(u)$ is absolutely convergent ($g(u)K_\varepsilon(u)$ is bounded)\cr
\+&$G(u) = G(u+) + G(u-)$\cr
\+&$G$ has {\it locally} bounded variation (as $V_a^\infty(g) = O(a^{-1-\varepsilon})$ uniformly in $a \geq a_0 > 0$)\cr
\+&$\sum_{p,k} \log(p)\,\int_0^\infty u^c |g(p^k u)| \, {du\over u} < \infty$ ($1<c<1+\varepsilon)$\cr
so that $\widehat{G}(s)$ exists for $1 < Re(s) <1+\varepsilon$ and that the conditions are met for the validity of Fourier-Mellin inversion on each such vertical line. The last item is then the most important
$$\widehat{G}(s) = -{\zeta'(s)\over\zeta(s)}\widehat{g}(s) \hbox{ for }1 < Re(s) <1+\varepsilon$$
and clinches the proof.

Let $L_r(s) = \pi^{-s/2} \Gamma(s/2)$, the question of the existence of
$$V_r(g) = \lim_{T\to\infty} {1\over{2\pi i}}\int_{c-iT}^{c+iT}-{L_r'(s)\over L_r(s)}\widehat{g}(s)\, ds$$
is more delicate (although we can now allow any $1+\varepsilon> c>0$)

One has (in this generality) to pair $g$ with its transpose $x \mapsto g^\tau(x) = {1\over x}g({1\over x})$ (so that $\widehat{g^\tau}(s) = \widehat{g}(1 - s)$ when both sides are defined) and consider
$$W_r(g) = \lim_{T\to\infty} {1\over{2\pi i}}\int_{c-iT}^{c+iT}-{L_r'(s)\over L_r(s)}(\widehat{g}(s)+\widehat{g}(1-s))\, ds$$
To study this we use the partial fraction expansion of $-{\Gamma'(s)\over \Gamma(s)}$ which gives
$$\eqalign{
-{L_r'(s)\over L_r(s)} &= (\log(\pi) + \gamma)/2 + {1\over s} - s\sum_{j\geq1}{1\over 2j(s+2j)} \cr
-{L_r'(s)\over L_r(s)} &= (\log(\pi) + \gamma)/2 + \int_1^\infty t^{-s}\, {dt\over t} + s\int_1^\infty {1\over 2}\ln(1-t^{-2})t^{-s}\, {dt\over t} \cr 
}$$

Still assuming that $g$ satisfies conditions $(A_\varepsilon)$ and $(B)$, let
$$H(u) = \int_1^\infty (g(ut)+g^\tau(ut))\, {dt\over t}$$
$$K(u) = \int_1^\infty {1\over 2}\ln(1-t^{-2})(g(ut)+g^\tau(ut))\, {dt\over t}$$
\+\cr
\+{\bf Lemma:}&&&$H(u) = O(1), u \rightarrow 0$ \cr
\+&&&$H(u) = O(u^{-1-\varepsilon}), u \rightarrow \infty$\cr
\+&&&$K(u) = O(\sup(u^\varepsilon, u^2\ln(u))), u \rightarrow 0$\cr
\+&&&$K(u) = O(u^{-1-\varepsilon}), u \rightarrow \infty$\cr
The proof is left to the reader.

Fubini's theorem and $t^{-s}\widehat{g}(s) = \int_0^\infty u^sg(ut)\, {dt\over t}$ then give
$$-{L_r'(s)\over L_r(s)}(\widehat{g}(s)+\widehat{g}(1-s)) = {\log(\pi) + \gamma \over 2} (\widehat{g}(s)+\widehat{g}(1-s)) + \widehat{H}(s) + s\widehat{K}(s)$$

It is time to introduce condition $(C_c)$\par
$(C_c)\;K(u) = \int_u^\infty L(t)\, {dt\over t}$ for some $L(t)$ with $\int_0^{\infty} t^c|L(t)|\, {dt\over t} < \infty$, continuous at $1$ (so that $K'(1)$ exists and is $-L(1)$), and of bounded variation in a neighborhood of $1$.

{\bf Lemma:} Under condition $(C_c)$\quad$s\widehat{K}(s) = \widehat{L}(s)$ on $Re(s) = c$.

{\bf Proof:} $K(u) = \int_1^\infty L(ut)\, {dt\over t}$ so that $\widehat{K}(s) = \int_1^\infty (\int_0^\infty u^s L(ut)\, {du\over u})\, {dt\over t} = 
\int_1^\infty t^{-s} \widehat{L}(s)\, {dt\over t} = {1\over s}\widehat{L}(s)$, the integrals being manipulated according to Fubini's theorem.

{\bf Theorem:} Let $g$ be such that it satisfies $(A_\varepsilon), (B), (C_c)$. Then
$$\lim_{T\to\infty} {1\over{2\pi i}}\int_{c-iT}^{c+iT}-{L_r'(s)\over L_r(s)}(\widehat{g}(s)+\widehat{g}(1-s))\, ds$$
exists and its value $W_r(g)$ is $(\log(\pi)+\gamma)g(1) + H(1) - K'(1)$.

{\bf Proof:} The conditions are satisfied for Fourier-Mellin inversion of each term: for $\widehat{g}(s)$ and $\widehat{H}(s)$ thanks to $(A_\varepsilon), (B)$ ($H$ is continuous with left and right derivatives), for $s\widehat{K}(s)$ thanks to $(C_c)$.

The von\thinspace Mangoldt case is included:

{\bf Theorem:} The function $g(u) = 1$ if $1 < u < X$, $g(1) = g(X) = 1/2$, $g(u) = 0$ for other $u$'s satisfies $(A_\varepsilon), (B), (C_c)$ and one has $(X >1)$
$$W_r(g) = (\log(\pi) + \gamma)/2 + \ln(X) + {1\over 2}\ln(1 - X^{-2})$$

{\bf Proof:} One obtains the following values for $K(u) = \int_1^\infty {1\over 2}\ln(1-t^{-2})(g(ut)+g^\tau(ut))\, {dt\over t}$
\+ \cr
\+$u \geq X$&&&$K(u) = 0$ \cr
\+$1\leq u \leq X$&&&$K(u) = {1\over 2}\int_{u/X}^1 \ln(1 - t^2)\, {dt\over t}$ \cr
\+$1/X \leq u \leq 1$&&&$K(u) = {1\over 2}\int_{u/X}^u\ln(1-t^2)\, {dt\over t} + {1\over 2u}\int_u^1\ln(1-t^2)\, dt $\cr
\+$u \leq 1/X$&&&$K(u) = {1\over 2}\int_{u/X}^u\ln(1-t^2)\, {dt\over t} + {1\over 2u}\int_u^{uX}\ln(1-t^2)\, dt $\cr
from which one checks $-K'(1) = {1\over 2}\ln(1 - X^{-2})$.

It is possible to give $W_r(g)$ as a more explicit finite integral if $g$ satisfies some additional assumptions.

{\bf Theorem:} If $g$ admits a derivative, continous and locally of bounded variation, with $|g(u)| + |ug'(u)| = O(1/K_\varepsilon(u))$, then $W_r(g) = V_r(g) + V_r(g^\tau)$ with
$$\eqalign{
V_r(g) &= \lim_{T\to\infty} {1\over{2\pi i}}\int_{c-iT}^{c+iT}-{L_r'(s)\over L_r(s)}\widehat{g}(s)\, ds \cr
&= {(\log(\pi)+\gamma)\over 2}g(1) + \int_1^\infty g(t)\, {dt\over t} + \int_1^\infty {1\over t^2 - 1}(g(t) - g(1))\, {dt \over t}
}$$

{\bf Notes:}
(1) Weil's own formulation in [5] is in terms of a certain limit, and applies to a large class of $g$'s including the step-function necessary for von\thinspace Mangoldt's theorem. In the case of a $C^1$ function Weil does not require $g'(u)$ to have locally bounded variation.\hfil\break
(2) This finite form of Weil's term is given in Haran [8]. Another finite form was obtained by Barner [9]. There are infinitely many finite forms, depending on the chosen regularisation at $1$: a few more will be given later.

{\bf Proof:} One checks easily that $(A_{\varepsilon'}) (0<\varepsilon'<\varepsilon), (B)$ hold. For $(C_c)$ let's take one ``half" of $K(u)$
$$\eqalign{
M(u) &= \int_1^\infty {1\over 2}\ln(1-t^{-2})g(ut)\, {dt\over t} \cr
uM'(u) &= \int_1^\infty {1\over 2}\ln(1-t^{-2})ug'(ut)\, dt \cr
-uM'(u) &= \int_1^\infty {1\over t^2 - 1} (g(ut) - g(u))\, {dt\over t} \cr
}$$
The right-hand side $L(u)$ is continuous, $O(u^{-1-\varepsilon})$ for $u \rightarrow \infty$ and $O(\sup(u,u^\varepsilon\ln(u)))$ for $u \rightarrow 0$. To show that it is locally of bounded variation one uses a third representation:
$$L(u) = - \int_1^2 {1\over 2}\ln(1-t^{-2}) utg'(ut)\, {dt\over t}
+ {1\over 2}\ln({3\over4})(g(2u) - g(u))
+ \int_2^\infty {1\over t^2 - 1} (g(ut) - g(u))\, {dt\over t}$$
So all conditions are met and the proof is complete.

{\bf Theorem (``Explicit Formula''):} Let $g$ be such that $(A_\varepsilon), (B)$ holds. Then, for any $0<c<1+\varepsilon$
$$W_r(g) = \lim_{T\to\infty} {1\over{2\pi i}}\int_{c-iT}^{c+iT}-{L_r'(s)\over L_r(s)}(\widehat{g}(s)+\widehat{g}(1-s))\, ds$$
exists if and only if 
$$Z(g) = \widehat{g}(0) + \widehat{g}(1) - \lim_{T\to\infty}\sum_{|Im(\rho)| < T}{\widehat{g}(\rho)}$$
exists and in case they do one has
$$Z(g) = \sum_p (V_p(g) + V_p(g^\tau)) + W_r(g)$$

{\bf Proof:} The limit if it exists is independant of $c$ so we will take it $>1$. Let $Z(s) = L_r(s)\zeta(s)$ be the completed zeta function, so that $Z(1-s) = Z(s)$. We need the following classical estimates (see [4], chap. 3):\hfil\break
(1) ${Z'(s)\over Z(s)} = O(\log(|Im(s)|))$ on $Re(s) = 1-c, Re(s) = c, |Im(s)| >2$\hfil\break
(2) There are $O(\log(T))$ zeros with $T\leq|Im(\rho)|\leq T+1$\hfil\break
(3) If $T$ is not the imaginary part of a zero then on ${-1\over 2}\leq Re(s)\leq {3\over 2}\; , |Im(s)| = T$
$$|{Z'(s)\over Z(s)}| = O(\log(T)) + O\left(\log(T)(\sup\left|{1\over T - |Im(\rho)|}\right|)\right)$$
so that there is a sequence $T_n$ going to infinity with $T_{n+1}\leq T_n + 1$ and $|{Z'(s)\over Z(s)}| = O(\log(T_n)^2)$
on ${-1\over 2}\leq Re(s)\leq {3\over 2}\; , |Im(s)| = T_n$.

As $|\widehat{g}(s)| = O(1/|s|)$ in  $1-c \leq Re(s) \leq c$, $W_r(g)$ will exist if and only if the limit through the $T_n$'s exists. Using the obvious contour integral and the functional equation one gets
$${1\over{2\pi i}}\int_{c-iT_n}^{c+iT_n}-{Z'(s)\over Z(s)}(\widehat{g}(s)+\widehat{g}(1-s))\, ds = \widehat{g}(0) + \widehat{g}(1) - \sum_{|Im(\rho)| < T_n}{\widehat{g}(\rho)} + O(\log(T_n)^2/T_n)$$
The theorem follows.

Applying this to the case of $\widehat{g}(s) = (X^s - 1)/s$ one gets on one side
$$\sum W_\nu(g) = \sum_{1<n<X}{\Lambda(n)} + {1\over2}\Lambda(X) + (\log(\pi) + \gamma)/2 + \log(X) + {1\over 2}\log(1 - X^{-2})$$
and on the other side
$$\log(X) + X - 1 - \lim_{T\to\infty}\sum_{|Im(\rho)| < T}{(X^\rho - 1)\over\rho}$$
so that
$$\sum_{1<n<X}{\Lambda(n)} + {1\over2}\Lambda(X) = X - \lim_{T\to\infty}\sum_{|Im(\rho)| < T}{(X^\rho - 1)\over\rho} - 1 - (\log(\pi) + \gamma)/2 - {1\over 2}\log(1 - X^{-2})$$
This reduces to von\thinspace Mangoldt's formula if one uses
$$\sum_{\rho}{1\over{\rho(1-\rho)}} = 2 + \gamma - \log(4\pi)$$

{\bf A convolution algebra}

The proof of the Explicit Formula would be quite easier if the complex integrals involved were absolutely convergent. This is possible only if $g$ is continuous, so it excludes the step-function necessary for von\thinspace Mangoldt's Theorem. Nevertheless Weil's positivity criterion tells us that it would be enough to understand what happens with smooth, compactly supported functions, so that we can restrict to continuous functions with not too many regrets.

Here is a class of such functions, stable under convolution, for which an Explicit Formula can be proven:

{\bf Theorem:} Let $(E)$ be the class of continuous functions $g$, such that $gK_\varepsilon$ is of bounded total variation on $(0, \infty)$ for some $\varepsilon > 0$, and such that $\int_{Re(s) = c} |\widehat{g}(s)| \log(|Im(s)|)\, ds < \infty$ for all $c$'s such that $-\varepsilon < c < 1 + \varepsilon$. Then $(E)$ is stable under multiplicative convolution and under the transpose $g \mapsto g^\tau(x) = {1\over x}g({1\over x})$. For each $g$ in $(E)$ the following explicit formula holds
$$\widehat{g}(0) + \widehat{g}(1) - \lim_{T\to\infty}\sum_{|Im(\rho)| < T}{\widehat{g}(\rho)} = \sum_\nu W_\nu(g)$$
with $$W_p(g) = \log(p)\sum_{k \geq 1} g(p^k) + \log(p)\sum_{k \geq 1} {1\over p^k}g({1\over p^k})$$
$$\eqalign{
W_r(g) = (\log(\pi) + \gamma)g(1) &+ \int_1^\infty (g(u) + g^\tau(u))\,{du \over u} \cr
+ \sum_{j\geq 1} \int_1^\infty (g(u) - g(1))u^{-2j}\,{du \over u}
&+ \sum_{j\geq 1} \int_1^\infty (g^\tau(u) - g(1))u^{-2j}\,{du \over u}
}$$
The sums defining the $W_\nu(g)$'s are absolutely convergent.

{\bf Proof:} Let $g$ and $k$ in $(E)$ and $\varepsilon$ working for both. To prove that $u^c\;(g*k)(u)$ has bounded total variation on $u \geq 1$ when $1 < c < 1 + \varepsilon$, we pick $c < d < 1 + \varepsilon$ and apply Mellin inversion on the line $Re(s) = d$
$$(g*k)(u) = {1\over 2\pi i}\int_{Re(s) = d} \widehat{g}(s) \widehat{k}(s)u^{-s}\, ds$$
Observing that $u \mapsto u^{c - d - iIm(s)}$ has a total variation on $u \geq 1$ which is $O(1 + |Im(s)|)$ and that $\left|\int_{Re(s) = d} \left|\widehat{g}(s) \widehat{k}(s) Im(s)\right|\, ds\right| < \infty$ completes the proof. Replacing $g$ and $k$ by $g^\tau$ and $k^\tau$ we get the same conclusion for $u^{-\varepsilon'} (g*k)(u)$ on $u\leq 1$ with $0 < \varepsilon' < \varepsilon$. The other conditions are trivially met. So $(E)$ is indeed stable under convolution.

For the existence of the explicit formula we just need to remark that $\int_{Re(s) = c}-{Z'(s)\over Z(s)}\widehat{g}(s)\, ds$ is absolutely convergent ($1 < c < 1 + \varepsilon$). And the value of $W_r(g)$ is obtained from
$$-{L_r'(s)\over L_r(s)} = (\log(\pi) + \gamma)/2 + {1\over s} + \sum_{j\geq1}({1\over s+2j} - {1\over 2j})$$
which can be used term by term thanks to the estimate on $Re(s) = c > 0$
$$\sum_{j\geq1}|{1\over s+2j} - {1\over 2j}| = O(\log(|s|))$$
so the proof will be complete if we check ${1\over 2\pi i}\int_{Re(s) = c} {\widehat{g}(s)\over s+2j}\, ds = \int_1^\infty g(u) u^{-2j}\, {du\over u}$\quad. This follows from interverting integrals starting with ${1\over s+2j} = \int_1^\infty u^{-s-2j}\, {du\over u}$, and applying Mellin inversion.

{\bf Weil's Theorem}

From now on we will take $g$ in the class $\cal D$ of smooth, compactly supported functions on $(0, \infty)$. Weil's Theorem is the simple, but pregnant observation that the contributions $W_\nu(g)$ to the Explicit Formula become truly identical in functional form when they are expressed as integrals on the completions of $\QQ$, the $\QQ_p's$ and $\RR$.

{\bf Theorem (Weil [5]):} For each completion $\QQ_\nu$ of $\QQ$ let $g_\nu$ be $x \mapsto g(|x|_\nu)$. Then
$$W_\nu(g) = PF_\nu \int_{\QQ_\nu^\times} {g_\nu({1\over x}) \over |1 - x|_\nu}\, d^\times x$$ 
where $d^\times x$ is the multiplicative Haar measure on $\QQ_\nu^\times$ which assigns $\log(p)$ as volume to the units if $\nu = p$ and is push-forwarded to $du\over u$ under $x \mapsto |x|$ for $\nu = r$ ($d^\times x = {dx\over 2|x|}$). The symbol $PF_\nu$ stands for a suitable ``finite part'' regularizing the singularity at $1$ of the integrand.

{\bf Proof:} For a finite prime we have (as $|1 - x|_\nu = 1$ for $|x|_\nu < 1$, $=|x|_\nu$ for $|x|_\nu > 1$)
$$W_p(g) = \log(p)\sum_{k \geq 1} g(p^k) + \log(p)\sum_{k \geq 1} {1\over p^k}g({1\over p^k})= \int_{|x|_\nu \neq 1} {g_\nu({1\over x}) \over |1 - x|_\nu}\, d^\times x $$
The manipulations are a little more magical for the real place (if $g(1) \neq 0$):
$$\eqalign{
W_r(g) &= V_r(g) + V_r(g^\tau) \cr
V_r(g^\tau) &= {\log(\pi) + \gamma \over 2}g(1) + \int_1^\infty {1\over u}g({1\over u})\,{du \over u} + \int_1^\infty {{1\over u}g({1\over u}) - g(1)\over u^2 - 1}\,{du \over u} \cr
V_r(g^\tau) &= {\log(\pi) + \gamma \over 2}g(1) + \int_1^\infty \left( {u\over u^2 -1}g({1\over u}) - {g(1) \over u^2 - 1} \right)\,{du \over u} \cr
V_r(g^\tau) &= {\log(\pi) + \gamma \over 2}g(1) + \int_1^\infty {g({1\over u}) - g(1) \over u - 1}\,{du \over 2u} + \int_1^\infty {g({1\over u})\over u + 1}\,{du \over 2u} + {g(1)\over 2}\int_1^\infty{du\over u(u+1)} \cr
V_r(g^\tau) &= {\log(2\pi) + \gamma \over 2}g_r(1) + \int_1^\infty {g_r({1\over x}) - g_r(1) \over |1 - x|}\,d^\times x + \int_{-\infty}^{-1} {g_r({1\over x}) \over |1 - x|}\,d^\times x \cr
V_r(g) &= {\log(2\pi) + \gamma \over 2}g_r(1) + \int_0^1 {g_r({1\over x}) - xg_r(1) \over |1 - x|}\,d^\times x + \int_{-1}^{0} {g_r({1\over x}) \over |1 - x|}\,d^\times x \cr
W_r(g) &= (\log(2\pi) + \gamma)g_r(1) + \int_{(0,2)} {g_r({1\over x}) - xg_r(1) \over |1 - x|}\,d^\times x + \int_{(-\infty,0)\,\cup\,(2,\infty)} {g_r({1\over x})  \over |1 - x|}\,d^\times x \cr
}$$
Replacing $g$ by $g^\tau$ and making the change of variable $x \rightarrow {1\over x}$ also gives
$$W_r(g) = (\log(2\pi) + \gamma)g_r(1) + \int_{x > {1\over 2}} {g_r({1\over x}) - g_r(1) \over |1 - x|}\,d^\times x + \int_{x<  {1\over 2}} {g_r({1\over x})  \over |1 - x|}\,d^\times x
$$
This completes the proof of Weil's theorem.

{\bf Notes:} (1) The result obtained at the archimedean place suggests to rewrite $W_p(g)$ as
$$W_p(g) = \int_{|x|_\nu = 1} {g_\nu({1\over x}) - g_\nu(1) \over |1 - x|_\nu}\, d^\times x + \int_{|x|_\nu \neq 1} {g_\nu({1\over x}) \over |1 - x|_\nu}\, d^\times x $$
and as shown by Weil [5] this is the correct form when dealing with a Dirichlet L-series (in which case $g_v(x) = \chi_\nu(x)g(|x|_\nu)$ for some local character $\chi_\nu(x)$)\ .\hfil\break
(2) Almost identical manipulations of $W_r(g)$ are in the paper of Haran [8]\ .\hfil\break
(3) Even more surprising is the case of a number field having a complex place: then the expression of $W_\CC(g)$ as an integral on the positive half-line is {\it a posteriori} seen to actually come from $\CC^\times$ through the fibration $z \rightarrow u = z\overline{z}$\ .

Weil's definition of the symbol ``PF'' at the real place is in terms of a certain limit, hence the comparison with finite primes is less satisfying. Haran has given a reformulation of the various $W_\nu$'s enabling to write them all in an {\it exactly} identical way, and this is our next topic.

{\bf The functional equation and homogeneous distributions}

To explain Haran's result I will first review some properties of homogeneous distributions on the local fields, and this will also give an opportunity to recall some aspects of Tate's Thesis [10] (see also Iwasawa [11]).

So let $\QQ_\nu$ be one of the completion of $\QQ$. Let $\chi$ be a (multiplicative) continuous character on $\QQ_\nu^\times$. A tempered distribution $F(t)$ is said to have homogeneity $\chi$ is $F(xt) = \chi(x)|x|^{-1} F(t)$ (that is if $\int F(t) \varphi(t)\, dt = \chi(x) \int F(t) \varphi(xt)\, dt$ for all test-functions $\varphi(t)$). For example the additive Haar measure has homogeneity $|x|$ whereas the Dirac $\delta$ function has homogeneity $1$. If $\nu$ is finite the Dirac $\delta$ is the only (tempered) distribution with support at the origin, whereas for the real place one also has its derivatives $\delta^{(k)}$ which have homogeneity $x^{-k}$. Note that it is only for even $k$ that this homogeneity is a power of $|x|$.

{\bf Theorem:} There exists a unique (up to a multiplicative constant) homogeneous distribution with homogeneity $\chi$.

{\bf Proof:} The reader is referred to the book of Gel'fand,Graev,Pyatetskii-Shapiro [12] and the other paper of Weil [13]. An example of this is that there is no {\it homogeneous} distribution on $\QQ_\nu$ which coincides with $1/|x|$ away from the origin as it would have the same homogeneity as the Dirac $\delta$-function. Of course there are some regularizations but they break the homogeneity. On the other hand there {\it is} a homogeneous distribution on $\RR$ which coincides with $1/x$ away from the origin (so that it has homogeneity $sign(x)$), and this is Cauchy's Principal Value.

For $Re(s) > 0$ the function $x \mapsto |x|^{s-1}$ is locally integrable hence defines a (tempered) distribution $\Delta_s$ whose homogeneity is $|x|^s$.

{\bf Theorem:} $\Delta_s$, as a function of s, has a meromorphic continuation to the entire complex plane, never vanishing, and with simple poles at those values of $s$ where the unique distribution with homogeneity $|x|^s$ is local (supported at the origin). 

{\bf Proof:} Again I refer to [12], [13].

So this means that the poles of $\Delta_s$ are the same as the poles of the local Euler factor of the Riemann Zeta function: the ${1\over 1 - p^{-s}}$'s and $\pi^{-s/2}\Gamma(s/2)$. This was pointed out by Weil [13]. This suggests very strongly to try to look at adeles in the hope to establish a link between homogeneous distributions and the zeros. This is just the well-known Iwasawa-Tate set-up, so we are following history in reverse here.

Of course on the adeles $\AA$ there are some difficulties with $|t|$ as it vanishes almost surely! (dropping the real component and restricting to the integral adeles $\prod\ZZ_p$ one can use the Borel-Cantelli lemma and $\sum{1\over p} = \infty$ to see that almost surely $|t|_p < 1$ for infinitely many $p$'s so that $|t| = \prod |t|_p = 0$). The trick is to define $\Delta_s(\varphi)$ as an integral over the {\it ideles} $\int_{\AA^\times} \varphi(x) |x|^s d^*x$ where $d^*x$ is a suitably normalized multiplicative Haar measure. But this works only for $Re(s) > 1$ where we can understand this integral as an infinite product. The test-function $\varphi$ is a finite linear combination of products $\prod_\nu \varphi_\nu$ with, for almost all $\nu$'s, $\varphi_\nu$ the characteristic function of the $p$-adic integers.

{\bf Theorem (Tate):} The homogeneous distribution $\Delta_s$ on the adeles $\AA$ defined as above for $Re(s) > 1$ has a meromorphic continuation to the entire complex plane and its zeros and poles are exactly those of the (completed) Riemann Zeta Function. The following functional equation is satisfied (${\cal F}$ = Fourier Transform)
$${\cal F}(\Delta_s) = \Delta_{1-s}$$
from which one can deduce the functional equation of the Zeta Function. In fact $Z(s)$ is obtained by evaluating $\Delta_s$ against a suitable test-function (independant of $s$). Furthermore these statements hold for a general algebraic number field and character of its idele classes.

{\bf Proof:} The reader is referred to Tate [10]. See also Iwasawa [11].

It is sometimes said that Tate's Thesis showed that the deeper meaning of the Functional Equation is that it reflects the validity of the Poisson Summation Formula in the adelic setting, but I prefer to see Poisson Formula as a tool to establish the fundamental equality ${\cal F}(\Delta_s) = \Delta_{1-s}$\enspace .

After this brave excursion in the realm of adeles and ideles, we can go back to the local situation on $\QQ_\nu$ and ask: what is the Fourier Transform of $\Delta_s$?

{\bf Theorem:} The Fourier Transform on $\QQ_\nu$ of $\Delta_s$ is $\Gamma_\nu(s)\Delta_{1-s}$ for a certain meromorphic function in the complex plane $\Gamma_\nu(s)$, analytic and non-vanishing in the critical strip.

{\bf Proof:} See Tate [10], Gel'fand [12], Weil [13].

The ``Tate-Gel'fand-Graev'' Gamma function $\Gamma_\nu(s)$ is simply the ratio of the $\nu$-adic components of $Z(s)$ and $Z(1-s)$
$$\Gamma_p(s) = {1 - p^{-1 + s}\over 1 - p^{-s}}$$
$$\Gamma_r(s) = \pi^{{1\over 2}-s} {\Gamma({s\over 2})\over\Gamma({1-s\over 2})}$$
One way of formulating the Functional Equation is through saying that it means that the ``adelic'' Gamma function is identically $1$, but this is a misleading statement because, as we will see later, its ``logarithmic derivative'' is basically the Explicit Formula.

Restricting to the critical strip where $\Delta_s$ and $\Delta_{1-s}$ are given by bona fide locally integrable functions, the Theorem means that for an arbitrary test-function $\varphi(y)$ on $\QQ_\nu$ (in the Schwartz class for $\RR$, locally constant with compact support for $\QQ_p$), with Fourier transform $\widetilde{\varphi}(x)$ the following identity of analytic functions of $s$ holds:
$$\int \widetilde{\varphi}(x) |x|^{s-1} dx = \Gamma_\nu(s) \int \varphi(y) |y|^{-s} dy$$
{\bf Lemma:} The identity holds if one only assumes that both integrals are absolutely convergent (both $\varphi(y)$ and $\widetilde{\varphi}(x)$ are supposed to be measurable locally integrable functions, $\widetilde{\varphi}(x)$ is the Fourier transform of $\varphi(y)$ in the sense of distributions, and $0<Re(s)<1$).

{\bf Proof:} One checks that the change of variable of Tate's Thesis applies to the double integrals\hfil\break
\centerline{$\int \widetilde{\varphi}(x) |x|^{s-1} dx\cdot\int \omega(y) |y|^{-s} dy$ and $\int \varphi(y) |y|^{-s} dy\cdot\int \widetilde{\omega}(x) |x|^{s-1} dx$}
and hence shows that they are equal (for an arbitrary test-function $\omega$).

There is one last topic which pertains to this chapter, and this is the Fourier Transform $G$ of $-\log(|x|)$. It was determined over $\QQ_p$ by Vladimirov [14] and of course is a well-known result in the real case. One way to obtain it is to expand the fundamental equation ${\cal F}(\Delta_s) = \Gamma_\nu(s)\Delta_{1-s}$ around $s=1$ (once an explicit representation of $\Delta_s$ around $s = 0$ is available) and one gets (for the details see [6]):

{\bf Theorem:}
$$G(\varphi) = {\log(p) \over 1 - 1/p}\left(\int_{|t|\leq 1}(\varphi(t) - \varphi(0))\, {dt\over|t|} + \int_{|t|>1}\varphi(t)\, {dt\over|t|} + {1\over p}\varphi(0)\right)\leqno (\QQ_p)$$
$$G(\varphi) = \int_{|t| \leq 1}(\varphi(t) - \varphi(0))\, {dt \over 2|t|} + \int_{|t| > 1}\varphi(t)\, {dt \over 2|t|} +\ (\log(2\pi) + \gamma)\cdot\varphi(0)\leqno (\RR)$$

{\bf Haran's Theorem}

This is the discovery by Haran [8] that although the local Weil term $W_\nu(g)$ starts its life as an inverse Mellin transform
$$W_\nu(g) = {1\over{2\pi i}}\int_{Re(s) = c}-{L_\nu'(s)\over L_\nu(s)}(\widehat{g}(s)+\widehat{g}(1-s))\, ds$$
and hence is first obtained as a {\it multiplicative} convolution it also can be shown to be an {\it additive} convolution once it has been lifted to be $\nu-$adics as advocated by Weil.

Let $R_s$ be the kernel $\Delta_{s}\over\Gamma_\nu(s)$ on $\QQ_\nu$. It is holomorphic in $Re(s) < 1$ as $\Delta_s$ has the same poles there as $\Gamma_\nu(s)$. Its Fourier Transform is $|x|^{-s}$, so $\left.{\partial\over\partial s}\right|_{s=0}R_s$ is just the Fourier transform $G$ of $-\log(|x|)$. 

{\bf Theorem (Haran [8]):} Let $*$ be the symbol of additive convolution on $\QQ_\nu$.
$$W_\nu(g) = \left.{\partial\over\partial s}\right|_{s=0}(R_s*g_\nu)(1)$$
{\bf Proof:} One writes the right-hand side as $(G*g_\nu)(1) = G(g_\nu(1-t)) = G(g(|1-t|_\nu))$ and it is then a matter of using the formula for $G$ given in the previous section and comparing the outcome with the previously obtained formula for $W_\nu(g)$ as an integral over the $\nu-$adics.

Of course, this is a little miraculous, and I will now turn to the ideas relevant to a more direct proof ([6]).

{\bf Mellin and Fourier}

The Mellin inversion formula has been in constant use in this paper
$$\hbox{for }u\in(0,\infty)\quad{1\over{2\pi i}}\int_{Re(s) = c}\widehat{g}(s)u^{-s}\, ds = g(u)$$
There is a cousin to this identity which is not so well-known. As before starting with $g$ a smooth compactly supported function on $(0,\infty)$ we can define $g_\nu$ on $\QQ_\nu$ as $y \mapsto g(|y|_\nu)$. This function will be a Schwartz-function hence its Fourier transform $\widetilde{g_\nu}(x)$ will also be of Schwartz class (and will depend only on $|x|_\nu$).

{\bf Theorem:} for $0<c<1$, $x \neq 0$
$${\cal F}^{-1}({g_\nu})(x) = {1\over{2\pi i}}\int_{Re(s) = c}\widehat{g}(s){|x|_\nu^{s-1}\over\Gamma_\nu(s)}\, ds$$

{\bf Note:} Of course ${\cal F}^{-1}({g_\nu}) = \widetilde{g_\nu}$ as $g_\nu$ is ``even'', but this is the statement that generalizes when using the $\Gamma$ function of a non-trivial character (in that case $g_\nu$ is twisted by the character).

{\bf Proof:} The integral converges because $\widehat{g}(s)$ decreases faster than any inverse power of $s$ while $\Gamma_\nu(s)$ is periodic for $\nu = p$ and satisfies the estimate ${\Gamma_r(s)}^{-1} = \Gamma_r(1-s) = O(|s|^{1/2})$ for $\nu = r$ (an interesting exercise using Stirling's Formula). Both sides of the equation define continuous locally integrable functions of $x$ (the right-hand side is $O(|x|^{c-1})$, $x \neq 0$). It is enough to prove that they are equal as distributions:
$$\eqalign{
\int \left({1\over{2\pi i}}\int_{Re(s) = c}\widehat{g}(s){|x|_\nu^{s-1}\over\Gamma_\nu(s)}\,ds\right) \widetilde{\varphi}(x)\,dx &= {1\over{2\pi i}}\int_{Re(s) = c}\widehat{g}(s)\left(\int{|x|_\nu^{s-1}\over\Gamma_\nu(s)}\widetilde{\varphi}(x)\,dx\right)\,ds \cr
= {1\over{2\pi i}}\int_{Re(s) = c}\widehat{g}(s)\left(\int |y|_\nu^{-s}\varphi(y)\,dy\right)\,ds &=\int \left({1\over{2\pi i}}\int_{Re(s) = c}\widehat{g}(s)|y|_\nu^{-s}\,ds\right) \varphi(y)\,dy \cr
= \int g(|y_\nu|) \varphi(y)\,dy &= \int {\cal F}^{-1}(g_\nu)(x) \widetilde{\varphi}(x)\,dx
}$$

This means that one can compute a $p-$adic Fourier transform with a complex integral! In the real case it also establishes a connection between the Fourier transforms of $t\mapsto g(e^t)$ and of $x\mapsto g(|x|)$.

{\bf The explicit formula: adelic approach}

The title of this chapter is more ambitious that what has been really accomplished. What I will explain now is a method to prove the Explicit Formula which treats all the places in exactly the same way, and which directly gives a $\nu-$adic expression to $W_\nu$, but it is not quite adelic (yet).

Starting from
$$W_\nu(g) = {1\over{2\pi i}}\int_{Re(s) = 3/2}-{L_\nu'(s)\over L_\nu(s)}(\widehat{g}(s)+\widehat{g}(1-s))\, ds$$
one shifts the term containing $\widehat{g}(s)$ to the line $Re(s) = c$ and the term containing $\widehat{g}(1-s)$ to the line $Re(s) = 1-c$, with $0<c<1$ and then makes the change of variable $s \rightarrow 1 -s$ in that last integral to obtain:
$$W_\nu(g) = {1\over{2\pi i}}\int_{Re(s) = c}-{\Gamma_\nu'(s)\over \Gamma_\nu(s)}\widehat{g}(s)\, ds$$
So the Explicit Formula is the inverse Mellin transform of the logarithmic derivative of the ``adelic'' Gamma function. Let's write $\Lambda_\nu(s)$ for $-{\Gamma_\nu'(s)\over \Gamma_\nu(s)}$.

One then takes the derivative of the fundamental identity of integrals on $\QQ_\nu$ for an arbitrary test-function $\varphi(y)$:
$$\int \widetilde{\varphi}(x) |x|_\nu^{s-1}\,dx = \Gamma_\nu(s) \int \varphi(y) |y|_\nu^{-s}\,dy$$
and rearranging terms, one obtains ($0<Re(s)<1$)
$$\Lambda_\nu(s)\int \varphi(y) |y|_\nu^{-s}\,dy = \int \varphi(y) (-\log(|y|_\nu))|y|_\nu^{-s}\,dy
+ \int \widetilde{\varphi}(x) (-\log(|x|_\nu)){|x|_\nu^{s-1}\over\Gamma_\nu(s)}\,dx$$
Integrating both sides against ${1\over{2\pi i}}\int_{Re(s) = c}\widehat{g}(s)\cdot\ \,ds$ gives
$$\int\varphi(y) W_\nu(g;y)\,dy =\int \varphi(y) (-\log(|y|_\nu)g_\nu(y))\,dy + \int\widetilde{\varphi}(x) (-\log(|x|_\nu){\cal F}^{-1}({g_\nu})(x))\,dx$$
with
$W_\nu(g;y) = {1\over{2\pi i}}\int_{Re(s) = c}\Lambda_\nu(s)\widehat{g}(s)|y|_\nu^{-s}\, ds$ which is a continuous $O(|y|^{-c})$ function on $\QQ_\nu$, $y\neq 0$ (hence locally integrable), so that the previous equation which is valid for an arbitrary test-function $\varphi(y)$ can be converted into the pointwise identity ($y\neq0$)
$$W_\nu(g;y) = -\log(|y|_\nu)g_\nu(y) + {\cal F}\left(-\log(|x|_\nu)\cdot{\cal F}^{-1}({g_\nu})\right)(y)$$
(we know that the last term is continuous because it is the Fourier transform of an $L^1$-function)
$$W_\nu(g;y) = -\log(|y|_\nu)g_\nu(y) + (G_\nu*g_\nu)(y)$$
with $G_\nu$ the Fourier transform of $-\log(|x|_\nu)$ and $*$ the symbol of {\it additive} convolution on $\QQ_\nu$. At $y=1$ we end up with Haran's Theorem, hence we have obtained the local contribution to the Explicit Formula directly as a $\nu-$adic object, in a manner exactly identical for all places. Furthermore this method works for all Dirichlet-Hecke L-series ([6]).

{\bf The conductor operator}

Let $H$ be the operator from ${\cal S}(\QQ_\nu)$ to $L^2(\QQ_\nu, dt)$
$$H(\varphi)(t) = \log(|t|_\nu)\varphi(t) + {\cal F}\left(\log(|\xi|)\cdot{\cal F}^{-1}(\varphi)(\xi)\right)(t)$$
By construction $H$ commutes with the Fourier transform. I have called it the ``conductor operator'' for reasons having to do with ramified characters at finite places [6]. A simple computation, left to the reader, shows it has a very important property:

{\bf Theorem:} The conductor operator commutes with the isometric action of the multiplicative group $\varphi(t) \mapsto |u|^{1/2} \varphi(ut)$.

In other words $H$ is dilation invariant. At the real place a classical cousin to $H$ is the Hardy Transform. The Hardy Transform is invariant under positive dilations only, this is an essential difference. The complete dilation invariance of $H$ means that it can be analyzed in terms of multiplicative characters, and thus exhibited as a multiplicative convolution. Of course this is not news to us, we already know a presentation of $H$ as a multiplicative convolution: the Explicit Formula! Indeed it can be shown [15] that the Explicit Formula, when put on the {\it critical} line, realizes the spectral analysis of $H$ (for this one needs all the explicit formulae associated to all multiplicative characters, not only the trivial one).

On $\QQ_p$, restricting $H$ to the ``cuspidal'' subspace $L_0$ of $L^2$ (the kernel of averaging over the units $|u| = 1$) one obtains a surprising positivity property:

{\bf Theorem:} The ``cuspidal'' spectrum of $H$ is $\{\log(p), 2\log(p), 3\log(p) \dots\}$ for $p>2$ and $\{2\log(2), 3\log(2) \dots\}$ for $p=2$.

The proof is in [6]. The ``invariant'' spectrum of $H$ is continuous and bounded [15].

A further surprising property of the conductor operator, valid at every place, is that it commutes with the inversion: let ${\cal D}_0(\QQ_\nu)$ be the space of Schwartz-Bruhat functions with compact support  away from the origin. $H$ can be seen as a map ${\cal D}_0(\QQ_\nu) \rightarrow L^2(\QQ_\nu)$. The inversion $I: \varphi \mapsto (t \mapsto {1\over |t|}\varphi({1\over t}))$ acts on the domain and (isometrically) on the target.

{\bf Theorem:} $H$ commutes with the inversion $I$.

{\bf Proof:} I will prove it here for $\QQ_\nu = \RR$ and even functions only, the general case is done in [15]. The proof is just a matter of observing that with our previous notation $\varphi(t) = g(|y|), y = t$, one simply has $$H(\varphi)(t) = - W_r(g; y) = -{1\over{2\pi i}}\int_{Re(s) = c}\Lambda_r(s)\widehat{g}(s)|y|_\nu^{-s}\, ds$$ and this is also according to $\Lambda_r(s) = \Lambda_r(1-s)$ $$-{1\over{2\pi i}}\int_{Re(s) = 1-c}\Lambda_r(s)\widehat{g}(1-s)|y|_\nu^{s-1}\, ds = -{1\over |y|}W_r(g^\tau; {1\over y}) = I(H(I(\varphi)))(t)$$

So the conductor operator, which gives an operator theoretic and spectral interpretation of the Explicit Formulae of Analytic Number Theory, has at every place some remnants of a conformal invariance, restricted to preserve $\{0\}\cup\{\infty\}$.

{\bf The symmetries of the explicit formula}

The representation of Weil's local terms as an additive convolution with the Fourier transform of $-\log(|x|)$ (Haran's Theorem) seems at first sight to have a rather convincing rigidity. But at the same time it reveals the amount of hidden flexibility there was in Weil's formulae, as all Fourier transforms depend on the choice of a basic additive character. Of course we can not make completely arbitrary choices if we do want to get the sum over the zeros, and this restriction is the only truly adelic aspect of our result which is immediately apparent: the local additive characters have to be chosen so that globally the group of principal adeles remains its own annihilator in the adele ring for the duality corresponding to the ensuing global additive character. The set of global additive characters satisfying this condition is a torsor under the multiplicative group of the rational number field, which thus appears as a simple but non trivial symmetry of the Explicit Formula. Concretely this means that after the choice of a non-zero rational number $q$ the shifted local terms $W_\nu(g) + \log(|q|_\nu)g(1)$ are equally valid. This can also be understood if one realizes that there is nothing sacred about the decomposition in Euler factors, each of which can be modified by a multiplicative factor $|q|_\nu^{s}$ without affecting the Zeta Function. This looks odd for the Riemann Zeta Function, but appears more natural when one considers the Euler Factors for Dirichlet $L$-series in case of ramification. And there is also the interesting inversion symmetry which acts trivially on each local term, but this is an ``accident'' due to the fact that ${1\over 1} = 1$.

Goldfeld has given an interpretation of Weil's explicit formula which has as its basis the group generated by the inversion and the non-zero rational numbers [16]. Connes [17] has also met the need for being very specific on what exactly are Weil's local terms. I must refer the reader to these papers for further information on their author's methods.

My own conclusion ([6]) is that any attempt to solve the Riemann Hypothesis will have to incorporate, if not to put at its forefront, the dilation and inversion invariance examplified by the conductor operator and partially broken by the Explicit Formula. 

{
{\bf REFERENCES}\par
\baselineskip = 12 pt
\parskip = 4 pt
\font\smallRoman = cmr8
\smallRoman
\font\smallBold = cmbx8
\font\smallSlanted = cmsl8
{\smallBold [1] B. Riemann},{\smallSlanted ``\"Uber die Anzahl der Primzahlen unter einer gegebenen Gr\"osse''}, Monatsber. Akadem. Berlin, 671-680, (1859).\par
{\smallBold [2] H.M. Edwards}, {\smallSlanted ``Riemann's Zeta Function''}, Academic Press, (1974).\par
{\smallBold [3] A.E. Ingham}, {\smallSlanted ``The Distribution of Prime Numbers''}, Cambridge Univ. Press, (1932).\par
{\smallBold [4] S.J. Patterson}, {\smallSlanted ``An introduction to the theory of the Riemann Zeta-Function''}, Cambridge Univ. Press, (1988).\par
{\smallBold [5] A. Weil},{\smallSlanted ``Sur les ``formules explicites'' de la th\'eorie des nombres premiers''}, Comm. Lund (vol d\'edi\'e \`a Marcel Riesz), (1952).\par
{\smallBold [6] J.F. Burnol}, {\smallSlanted ``The Explicit Formula and a Propagator''}, electronic manuscript, available at the http://xxx.lanl.gov server, math/9809119 (September 1998, revised November 1998)\par
{\smallBold [7] Math.Soc.Japan}, {\smallSlanted ``Stationary Processes''}, Article 395 from the Encyclopedic Dictionary of Mathematics, 3rd ed., translation MIT Press (1987).\par
{\smallBold [8] S. Haran},{\smallSlanted ``Riesz potentials and explicit sums in arithmetic''}, Invent. Math.  101, 697-703 (1990).\par
{\smallBold [9] K. Barner}, {\smallSlanted ``On A.Weil's explicit formula''}, Journal f\"ur Mathematik Band 323, 139-152, (198?).\par
{\smallBold [10] J. Tate}, Thesis, Princeton 1950, reprinted in Algebraic Number Theory, ed. J.W.S. Cassels and A. Fr\"ohlich, Academic Press, (1967).\par
{\smallBold [11] K. Iwasawa}, {\smallSlanted ``Letter to J. Dieudonn\'e'' (1952)}, in ``Zeta Functions in Geometry'', Adv.Stud.Pure Math 21 ed. Kurokawa, Sunada, pub. Kinokuniya, (1992).\par
{\smallBold [12] I. M. Gel'fand, M. I. Graev, I. I. Piateskii-Shapiro},{\smallSlanted ``Representation Theory and automorphic functions''}, Philadelphia, Saunders (1969).\par
{\smallBold [13] A. Weil},{\smallSlanted ``Fonctions z\^etas et distributions''}, S\'eminaire Bourbaki n${}^{\smallRoman o}$ 312, (1966).\par
{\smallBold [14] V. S. Vladimirov},{\smallSlanted ``Generalized functions over the fields of p-adic numbers''}, Russian Math. Surveys 43:5, 19-64 (1988).\par
{\smallBold [15] J.F. Burnol}, {\smallSlanted ``Spectral analysis of the local conductor operator''}, math/9811040 (November 1998).\par
{\smallBold [16] D. Goldfeld}, {\smallSlanted ``A spectral interpretation of Weil's explicit formula''}, Lect. Notes Math 1593, p 135-152, Springer Verlag (1994).\par
{\smallBold [17] A. Connes}, {\smallSlanted ``Trace formula in non-commutative Geometry and the zeros of the Riemann zeta function''}, math/9811068 (November 1998).\par
\vfill
\centerline{Jean-Fran\c{c}ois Burnol, 62 rue Albert Joly, F-78000 Versailles, France}
\centerline{jf.burnol@dial.oleane.com}
\centerline{October 1998, revised November 1998}
}
\eject
\bye